\documentstyle[12pt,fleqn]{article}

\font\twlgot =eufm10 scaled \magstep1
\font\egtgot =eufm8
\font\sevgot =eufm7

\font\twlmsb =msbm10 scaled \magstep1
\font\egtmsb =msbm8
\font\sevmsb =msbm7

\newfam\gotfam

\textfont\gotfam\twlgot
\scriptfont\gotfam\egtgot
\scriptscriptfont\gotfam\sevgot

\newfam\msbfam
\textfont\msbfam\twlmsb
\scriptfont\msbfam\egtmsb
\scriptscriptfont\msbfam\sevmsb

\def\pBbb{\relax\ifmmode\expandafter\Bb\else\typeout{You cann't use
Bbb in text mode}\fi}
\def\Bb #1{{\fam\msbfam\relax#1}}

\def\thebibliography#1{\section*{References}\list
  {[\arabic{enumi}]}{\settowidth\labelwidth{#1}\leftmargin\labelwidth
    \advance\leftmargin\labelsep
    \usecounter{enumi}}
    \def\newblock{\hskip .11em plus .33em minus .07em}
    \sloppy\clubpenalty4000\widowpenalty4000
    \sfcode`\.=1000\relax}

\let\Large=\large

\def\op#1{\mathop{\fam0 #1}\limits}

\newcommand{\id}{{\rm Id\,}}

\newcommand{\di}{{\rm dim\,}}

\newcommand{\beq}{\begin{equation}}
\newcommand{\eeq}{\end{equation}}
\newcommand{\ben}{\begin{eqnarray}}
\newcommand{\een}{\end{eqnarray}}
\newcommand{\be}{\begin{eqnarray*}}
\newcommand{\ee}{\end{eqnarray*}}
\newcommand{\bea}{\begin{eqalph}}
\newcommand{\eea}{\end{eqalph}}

\newcommand{\cF}{{\cal F}}

\newcommand{\fl}{\flat}
\newcommand{\sh}{\sharp}

\newcommand{\al}{\alpha}

\newcommand{\Om}{\Omega}

\newcommand{\ar}{\op\longrightarrow}

\newcounter{eqalph}
\newcounter{equationa}
\newcounter{theorem}
\newcounter{proposition}
\newcounter{lemma}
\newcounter{corollary}
\newcounter{definition}

\def\thedefinition{\arabic{definition}}

\newenvironment{proof}{\noindent 
{\bf Proof.}}{\hfill $\Box$ \medskip}

\newenvironment{theo}{\refstepcounter{definition} \medskip\noindent
{\bf Theorem \thedefinition.}\it}{\medskip }

\newenvironment{eqalph}{\stepcounter{equation}
\setcounter{equationa}{\value{equation}}
\setcounter{equation}{0}

\begin{eqnarray}}{\end{eqnarray}\setcounter{equation}{\value{equationa}}}

\newcommand{\mar}[1]{}

\begin{document}
\hbox{}

{\parindent=0pt 

{\Large \bf  Recursion operators between degenerate Poisson structures}
\bigskip

{\sc Gennadi Sardanashvily}

{ \small 

{\it Department of Theoretical Physics,
Physics Faculty, Moscow State University, 117234 Moscow, Russia, 

E-mail:
sard@grav.phys.msu.su; sard@campus.unicam.it

URL: http://webcenter.ru/$\sim$sardan/}
\bigskip

{\bf Abstract.}

Two degenerate Poisson structures of the same rank possess a recursion
operator if and only if their characteristic distributions coincide.

\medskip

\bigskip
\bigskip

} }

Let $Z$ be a smooth real manifold provided with Poisson bivectors
$w$ and $w'$ of constant rank $k\leq\di Z/2$. They define the
bundle homomorphisms 
\mar{p4}\beq
w'^\sh: T^*Z\to TZ, \qquad w^\sh: T^*Z\to TZ \label{p4}
\eeq
over $Z$.
Let $R$ be a
tangent-valued one-form (a type (1,1) tensor field) on $Z$.
It yields bundle endomorphisms
\be
R: TZ\to TZ, \qquad R^*: T^*Z\to T^*Z.
\ee 
It is called a recursion operator if
\mar{p0}\beq
w'^\sh=R\circ w^\sh=w^\sh\circ R^*. \label{p0}
\eeq
The well-known sufficient conditions for $R$ to be a recursion 
operator (\ref{p0}) are the following: (i) the Nijenhuis torsion of $R$
vanishes, and (ii) the Magri--Morosi concomitant of $R$ and $w$ do so
\cite{1,2}. However, these conditions fail to characterize straightforwardly
the Poisson structures related by a recursion operator.

\begin{theo} \label{p1} \mar{p1}
A recursion operator (\ref{p0}) between regular Poisson structures
of the same rank exists if and only if their characteristic distribution
coincide.
\end{theo}

\begin{proof}
It follows from the equalities
(\ref{p0}) that a recursion operator $R$
sends the characteristic distribution of $w$ to 
that of $w'$, and these distributions coincide if
$w$ and $w'$ are of the same rank.

Let now $w$ and $w'$ be regular Poisson structures whose 
characteristic distributions
coincide. Let $\cF$ be their characteristic foliation 
and $i_\cF:T\cF\to TZ$ their characteristic
distribution. Given the dual $T\cF^*\to Z$ of $T\cF\to Z$,
we have the exact sequences of vector bundles
\mar{p2,3}\ben
&& 0\to T\cF \ar^{i_\cF} TZ \ar TZ/T\cF\to 0, \label{p2} \\
&& 0\to {\rm Ann}\,T\cF\ar T^*Z\ar^{i^*_\cF} T\cF^* \to 0 \label{p3}
\een
over $Z$.
Since Ann$\,T\cF\subset T^*Z$ is
the kernel of the bundle homomorphism $w^\sh$ and $w'^\sh$ (\ref{p4}), 
these factorize in a unique fashion 
\be
&& w^\sh:
T^*Z\ar^{i^*_\cF} T\cF^*\ar^{w^\sh_\cF}
T\cF\ar^{i_\cF} TZ, \\
&&  w'^\sh:
T^*Z\ar^{i^*_\cF} T\cF^*\ar^{w'^\sh_\cF}
T\cF\ar^{i_\cF} TZ 
\ee
through the bundle isomorphisms
\be
&& w_\cF^\sh: T\cF^*\to T\cF,  \qquad
w^\sh_\cF:\al\mapsto -w(z)\lfloor \al, \qquad
\al\in T_z\cF^*, \\
&& w'^\sh_\cF: T\cF^*\to T\cF,  \qquad
w'^\sh_\cF:\al\mapsto -w'(z)\lfloor \al, \qquad
\al\in T_z\cF^*. 
\ee
Let us consider the inverse isomorphisms 
\mar{p13}\beq
w_\cF^\fl : T\cF\to T\cF^*, \qquad w'^\fl_\cF : T\cF\to T\cF^* \label{p13}
\eeq
and the compositions
\mar{p10}\beq
R_\cF= w'^\sh_\cF\circ w_\cF^\fl: T\cF\to T\cF, \qquad
R_\cF^*= w_\cF^\fl \circ w'^\sh_\cF: T\cF^*\to T\cF^*. \label{p10}
\eeq
The relation
\mar{p11}\beq
w'^\sh_\cF=R_\cF\circ w^\sh_\cF=  w^\sh_\cF\circ R^*_\cF \label{p11}
\eeq
holds. In order to obtain a recursion operator (\ref{p0}), it suffices 
to extend the morphisms $R_\cF$ and $R_\cF^*$ (\ref{p10}) onto $TZ$ and $T^*Z$, respectively.
 For this purpose, let us consider a splitting
\be
\zeta: TZ\to T\cF, \qquad TZ=T\cF\oplus (\id-i_\cF\circ\zeta)TZ=T\cF\oplus E, 
\ee
of the exact sequence (\ref{p2}) and the dual splitting 
\be
\zeta^* :T\cF^*\to T^*Z, \qquad T^*Z=\zeta(T\cF^*)\oplus 
(\id-\zeta^*\circ i^*_\cF)T^*Z= \zeta^*(T\cF^*)\oplus E'
\ee
of the exact sequence (\ref{p3}). Then, the desired extensions are
\mar{p12}\beq
R:=R_\cF\times \id E, \qquad R^*:=(\zeta^*\circ R^*_\cF)\times \id E^*. \label{p12}
\eeq
The condition (\ref{p0}) is easily verified.
\end{proof}

The recursion operator $R$ (\ref{p12}) is by no means unique. It depends on the choice of  a splitting
of the exact sequences (\ref{p2}) -- (\ref{p3}), while its restriction $R_\cF$ (\ref{p10}) to
the characteristic distribution $T\cF$ is uniquely defined by the Poisson structures $w$ and $w'$.

Note that the inverse isomorphisms $w^\fl_\cF$ and $w'^\fl_\cF$ (\ref{p13}) yield the corresponding
leafwise symplectic forms $\Om_\cF$ and $\Om'_\cF$ on the foliated manifold $(Z,\cF)$ \cite{3,4}.
It enables one to generalize some results on recursion operators between symplectic structures
to Poisson ones in terms of the leafwise differential calculus.

\end{document}